\documentclass[11pt,reqno]{amsart}

\usepackage{multirow,multicol,soul,hhline,a4wide,amsmath,amssymb,dsfont,epsfig,graphicx,subfig,verbatim,multicol,enumitem,lipsum,amsaddr}
\usepackage[hang,flushmargin]{footmisc}

\usepackage[margin=2.5cm]{geometry}

\theoremstyle{definition}
\newtheorem{theorem}{Theorem}

\usepackage{cite}

\DeclareMathOperator{\argmin}{argmin}
\newcommand{\e}{\text{e}}
\newcommand{\DF}{\text{DF}}
\newcommand{\E}{\text{E}}

\newcommand{\cR}{\mathcal{R}}
\newcommand{\cS}{\mathcal{S}}
\newcommand{\cA}{\mathcal{A}}
\newcommand{\cB}{\mathcal{B}}
\newcommand{\cC}{\mathcal{C}}

\newcommand{\cP}{\mathcal{P}}
\newcommand{\cPnet}{\mathcal{P}^{\text{net}}}
\newcommand{\cPgross}{\mathcal{P}^{\text{gross}}}
\newcommand{\de}{\text{d}}
\newcommand{\dS}{\text{d}S}
\newcommand{\dI}{\text{d}I}
\newcommand{\dH}{\text{d}H}
\newcommand{\dD}{\text{d}D}
\newcommand{\dDt}{\text{d}D^\ast}
\newcommand{\dt}{\text{d}t}
\newcommand{\dN}{\text{d}N}
\newcommand{\tS}{\overline{S}}
\newcommand{\tI}{\overline{I}}
\newcommand{\tH}{\overline{H}}
\newcommand{\tD}{\overline{D}}

\newcommand{\ar}[1]{\multicolumn{1}{r|}{#1}}
\newcommand{\arr}[1]{\multicolumn{1}{rV{4}}{#1}}

\usepackage{pgf,tikz}
\usepackage{tkz-fct}
\usepackage{pgfplots}
\usetikzlibrary{arrows,trees}
\pgfplotsset{ every non boxed x axis/.append style={x axis line style=<->},
     every non boxed y axis/.append style={y axis line style=<->}}
\definecolor{newpurple}{RGB}{195, 22, 140}
\definecolor{newgray}{RGB}{240, 240, 240}
\definecolor{newlightblue}{RGB}{0, 175, 158}
\definecolor{newblue}{RGB}{47, 50, 145}
\definecolor{newyellow}{RGB}{232, 222, 0}
\definecolor{newgreen}{RGB}{0, 155, 1}
\usepackage{diagbox}
\usepackage{boldline}

\begin{document}
\title{Design and financial analysis of a health insurance based on an SIH-type epidemic model}
%\title{Premium computation and analysis of a health insurance based on an SIH-type epidemic model}
\author{\small Jonathan Hoseana, Felivia Kusnadi, Gracia Stephanie, Levana Loanardo, Catherine Wijaya}
\address{\normalfont\small Center for Mathematics and Society, Department of Mathematics, Parahyangan Catholic University, Bandung 40141, Indonesia}
\email{\,\,\,\begin{minipage}{0.9\textwidth}j.hoseana@unpar.ac.id, felivia@unpar.ac.id, hello.graciastephanie@gmail.com,\\ levanaloa03@gmail.com, catherinewij03@gmail.com\end{minipage}}
\date{}

\begin{abstract}
We present a design and financial analysis of a health insurance based on an SIH-type epidemic model. Specifically, we first construct the model in a continuous form, study its dynamical properties, and formulate the financial quantities involved in our insurance. Subsequently, we discretise the model using the forward Euler method, study the dynamical properties of the resulting discrete model, and formulate discrete analogues of the above financial quantities. We conduct a numerical simulation using two sets of parameter values, each representing a disease-free and an endemic scenario, which reveals that in the latter scenario, the insurance's gross premium is higher, the insurer's minimum loss-preventing start-up capital is lower, and the insurer's total profit is higher, compared to the corresponding values in the former scenario. Finally, through a sensitivity analysis, we show that in both scenarios, the disease's basic reproduction number, the gross premium, the minimum start-up capital, and the total profit depend most sensitively on the population's natural death coefficient, the disease's incidence coefficient, the hospitalisation benefit, and the premium surcharge percentage allocated to profit, respectively.


\end{abstract}

\maketitle

%-------------------------------------------------------------------------------------------------------------------------------

\section{Introduction}

In the wake of the recent coronavirus pandemic, the literature has witnessed a remarkable increase in studies on disease spread which utilise mathematical models. The majority of such studies, however, concentrate solely on modelling the dynamics of the spread itself, or on determining the successfulness of certain eradicative strategies. As a result, a major healthcare aspect which could also raise concern during disease outbreaks remains largely overlooked: the design and financial analysis of an accommodating health insurance.%As a result, a vital medical aspect which is also closely linked to epidemic spreading remains largely unexplored: the design of an appropriate health insurance. 

Indeed, studies involving the design of a health insurance during disease outbreaks have been relatively scarce. A pioneering work is that of Feng and Garrido \cite{Feng,FengGarrido} who, only in 2005, revisited the 1927 Kermack-McKendrick SIR-type epidemic model \cite{KermackMcKendrick} and devised a health insurance policy for individuals in the host population, applying the results to the 1665 great plague in England and to the 2003 SARS epidemic in Hong Kong. In 2020, this work was improved by Nkeki and Ekhaguere \cite{NkekiEkhaguere} through the construction of an SIDRS-type epidemic model, intended to provide a more comprehensive portrait of the disease's dynamics from an actuarial viewpoint. A notable advancement was the introduction of the deaths compartment, aimed to facilitate the actuarial formulations. Complementarily, Hainaut \cite{Hainaut}, in 2021, adopting a different approach, assumed that the time-evolution of the number of infected individuals obeys a specific functional form, which leads to the construction of a non-autonomous SIDS-type epidemic model characterised by a time-dependent incidence coefficient. An advantage of such an approach, as made apparent by Atatalab et al.\ \cite{AtatalabEtAl}, is that the majority of the calculated actuarial quantities admit explicit descriptions in terms of the gamma function \cite[p.\ 190]{Kreyszig}. Most recently, works along similar lines have been conducted by Nam \cite{Nam}, Guerra \cite{Guerra}, and Zhai et al.\ \cite{ZhaiEtAl}, employing, respectively, an SEIR-type, an SVIRD-type, and an SVEIRD-type epidemic model.

%Reviewing the works surveyed above, one readily discerns opportunities for further development. %In the works surveyed above, one readily identifies several limitations.

The works surveyed above possess several limitations. Firstly, none of the employed models allocate a compartment for hospitalised individuals, while in reality, health insurance policies often differentiate provisions for hospitalised infected individuals from those for non-hospitalised infected individuals, due to their different degrees of medical necessity. In addition, such policies may also distinguish the amount of benefit for a death caused by the disease under consideration, from that for a natural death. Furthermore, the calculations of the present values of premiums and benefits in the above works all utilised continuous interest-compounding, while in reality, interest-compounding generally occurs at discrete time-intervals, such as monthly. In fact, some of the works \cite{ZhaiEtAl} did not consider the effect of interest rates. Finally, each of the above works could be extended with a sensitivity analysis of the key actuarial quantities, such as the insurance's premium and the insurer's total profit, with respect to the parameters involved in the model, which could provide a quantitative assessment of the insurer's financial stability.

%Motivated by the above observations, in this paper we shall present a design and financial analysis of a health insurance based on an epidemic model.

The purpose of this paper is to present an epidemic-model-based design and financial analysis of a health insurance which resolves the above limitations. We organise our work as follows. Firstly, addressing the aforementioned absence of a hospitalised compartment, we construct an SIH-type continuous epidemic model, study its dynamical properties, and formulate the financial quantities involved in our insurance design based on the model (section \ref{sec:continuousmodel}). Subsequently, to prepare a setting for our numerical simulation, we discretise the model using the forward Euler method, study the discrete model's dynamical properties, and formulate the discrete analogues of the involved financial quantities (section \ref{sec:discretisation}). In our numerical simulation, we employ our discrete model to describe two scenarios characterised by two different sets of parameter values, each of which representing a disease-free and an endemic scenario (section \ref{section:numerical}). In each scenario, we calculate the insurance's monthly premium, the insurer's minimum loss-preventing start-up capital, and the insurer's end-of-period total profit. We also analyse the sensitivity of these quantities, and of the model's basic reproduction number, with respect to the existing parameters, in each of the two scenarios. Finally, we state our conclusions and describe avenues for future research (section \ref{sec:conclusions}). %Specifically, addressing ..., we shall first construct ...

\section{The model}\label{sec:continuousmodel}

In this section, we construct our continuous epidemic model, which we later use to design and analyse financially our health insurance. As mentioned in the previous section, our epidemic model is of type SIH, meaning that at any given time, each individual living in the population under consideration is assumed to possess exactly one of three different statuses: susceptible, infected, and hospitalised.

After detailing our model's construction (subsection \ref{subsec:construction}), we shall analyse the model from the perspective of dynamical systems theory \cite{Allen,Martcheva,Robinson}. Our aim is three-fold, namely, to establish the non-negativity and the boundedness of the model's solutions associated to non-negative initial conditions (subsection \ref{subsec:nonneg}), to determine the model's equilibria and basic reproduction number (subsection \ref{subsec:equilibria}), and to establish conditions for the local asymptotic stability of the model's equilibria, in connection to the model's basic reproduction number (subsection \ref{subsec:stability}). Finally, we design our health insurance based on the model (subsection \ref{subsec:design}).

\subsection{Model construction}\label{subsec:construction}

Suppose that a certain infectious disease spreads across a population. Suppose also that, at any given time, each individual living in the population can be assigned exactly one of the following three statuses: susceptible, infected, and hospitalised. Furthermore, let us denote by $S=S(t)$, $I=I(t)$, and $H=H(t)$ the number of susceptible, infected, and hospitalised individuals at time $t\geqslant 0$. To construct our model, we assume that the following occur with regards to the changes in these numbers at any given time.
\begin{enumerate}[leftmargin=0.9cm,itemsep=3pt]
\item[(i)] The number $S$ of susceptible individuals increases due to births occurring at the rate of $\lambda>0$, and the hospitalised and infected individuals recovering at the rates of $\alpha_1 H$ and $\alpha_2 I$, respectively, where $\alpha_1,\alpha_2>0$. The same number decreases due to infections occurring at the rate of $\beta S I$, where $\beta>0$, and natural deaths occurring at the rate of $\mu_1 S$, where $\mu_1>0$.

\item[(ii)] The number $I$ of infected individuals increases due to infections occurring at the rate of $\beta S I$, and decreases due to treatment occurring at the rate of $\alpha_2 I$, hospitalisation occurring at the rate of $\gamma I$, and deaths due to the disease occurring at the rate of $\mu_2 I$, where $\mu_2>0$.

\item[(iii)] The number $H$ of hospitalised individuals increases due to hospitalisation occurring at the rate of $\gamma I$. The same number decreases due to treatment occurring at the rate of $\alpha_1 H$, and deaths due to the disease occurring at the rate of $\mu_2 H$, where $\mu_2>0$.
\end{enumerate}
We remark that the choice of $\beta SI$ as the incidence rate assumes that hospitalised individuals are completely quarantined, implying no possibility of disease transmission from hospitalised individuals to susceptible individuals.

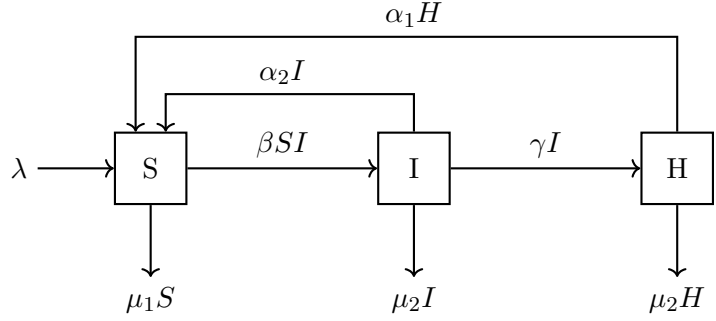
\begin{figure}
\centering
\begin{tikzpicture}
\node[rectangle,fill=white,draw=black,thick,minimum size=0.95cm] (S) at (0,0) {S};
\node[rectangle,fill=white,draw=black,thick,minimum size=0.95cm] (I) at (3.5,0) {I};
\node[rectangle,fill=white,draw=black,thick,minimum size=0.95cm] (H) at (7,0) {H};

\node (lambda) at (-1.75,0) {$\lambda$};
\draw[->,thick] (lambda) edge (S);

\node (muS) at (0,-1.75) {$\mu_1 S$};
\draw[<-,thick] (muS) edge (S);

\node (muI) at (3.5,-1.75) {$\mu_2 I$};
\draw[<-,thick] (muI) edge (I);

\node (muH) at (7,-1.75) {$\mu_2 H$};
\draw[<-,thick] (muH) edge (H);

\draw[->,thick] (S) edge node[above] {$\displaystyle\beta S I$} (I);

\draw[->,thick] (I) edge node[above] {$\gamma I$} (H);

\draw[->,thick] (H) -- (7,1.75) -- (-0.2,1.75) -- (-0.2,0.4875);
\node[above] at (3.5,1.75) {$\alpha_1 H$};

\draw[->,thick] (I) -- (3.5,0.9875) -- (0.2,0.9875) -- (0.2,0.4875);
\node[above] at (1.75,0.9875) {$\alpha_2 I$};
\end{tikzpicture}
\caption{\label{fig:compartment} The compartment diagram of our SIH-type model.}
\end{figure}

The above assumptions lead to the compartmental diagram in Figure \ref{fig:compartment}. We thus construct the SIH-type model
\begin{equation}\label{eq:model3eqns}
\left\{\begin{array}{rcl}
\displaystyle\frac{\dS}{\dt} &\!\!\!=\!\!\!& \displaystyle\lambda-\beta SI +\alpha_1 H +\alpha_2 I-\mu_1 S,\\[0.3cm]
\displaystyle\frac{\dI}{\dt} &\!\!\!=\!\!\!& \displaystyle\beta SI -\left(\alpha_2+\gamma+\mu_2\right)I,\\[0.3cm]
\displaystyle\frac{\dH}{\dt} &\!\!\!=\!\!\!& \displaystyle\gamma I-\left(\alpha_1+\mu_2\right)H.
\end{array}\right.
\end{equation}
The parameters involved in the model, along with their values used in our numerical simulations, are summarised in Table \ref{tab:parameters1}. Of note is that we use months as the unit for time in our model. In the upcoming subsections, we shall study the dynamical properties of the above model.

\begin{table}\centering
	\renewcommand{\arraystretch}{1.9}\renewcommand{\tabcolsep}{2pt}
		\scalebox{0.925}{\begin{tabular}{|c|c|c|c|c|} 
			\hline
			Parameter & Description & Unit & 
\begin{minipage}{1.9cm}\centering\smallskip
Value for simulation\medskip
\end{minipage} & Source \\
			\hhline{|=|=|=|=|=|}
			$\lambda$ & birth rate & individual/month & 4.21492 & \cite{UN,WHO,Affairs}\\\hline
			$\alpha_1$ & \begin{minipage}{6.8cm}\centering\smallskip treatment coefficient of hospitalised individuals\medskip\end{minipage} & 1/month & 0.05000 & \cite{BettgerEtAl}\\\hline
			$\alpha_2$ & \begin{minipage}{6.8cm}\centering\smallskip treatment coefficient of non-hospitalised infected individuals\medskip\end{minipage} & 1/month & 0.05000  & \cite{BettgerEtAl}\\\hline
			$\beta$ & incidence coefficient & 1/(individual\,$\times$\,month) & \begin{minipage}{1.9cm}\centering\smallskip 0.00100\\0.00300\smallskip\end{minipage} & simulated\\\hline
			$\gamma$ & \begin{minipage}{6.8cm}\centering\smallskip hospitalisation coefficient of infected individuals\medskip\end{minipage} & 1/month & 0.66000  & \cite{MalliaJohnston}\\\hline
			$\mu_1$ & natural death coefficient & 1/month & 0.00745 & \cite{SOA}\\\hline
			$\mu_2$ & death-by-disease coefficient & 1/month & 0.01829 & \cite{MalhotraEtAl}\\\hline
			$S(0)$ & initial number of susceptible individuals & individual & 2999 & simulated\\\hline
			$I(0)$ & initial number of infected individuals & individual & 1 & simulated\\\hline
			$H(0)$ & initial number of hospitalised individuals & individual & 0 & simulated\\\hline
		\end{tabular}}\smallskip
		\caption{\label{tab:parameters1}Epidemiological parameters used in our model and their values selected for our numerical simulations (section \ref{section:numerical}).}
\end{table}

\subsection{Non-negativity and boundedness of solutions}\label{subsec:nonneg}

\sloppy Consider a solution $\left(S,I,H\right)=\left(S(t),I(t),H(t)\right)$ of the model \eqref{eq:model3eqns}, associated to an initial condition $\left(S(0),I(0),H(0)\right)\in\mathbb{R}_+^3$, where $\mathbb{R}_+:=[0,\infty)$. By the model's first, second, and third equation, respectively, we have that, at every time $t\geqslant 0$ for which $S(t)=0$ we have
$$\frac{\dS}{\dt}=\lambda +\alpha_1 H +\alpha_2 I> 0,$$
at every time $t\geqslant 0$ for which $I(t)=0$ we have
$$\frac{\dI}{\dt}=0\geqslant 0,$$
and at every time $t\geqslant 0$ for which $H(t)=0$ we have
$$\frac{\dH}{\dt}=\gamma I\geqslant 0.$$
Since $\left(S(0),I(0),H(0)\right)\in\mathbb{R}_+^3$, these imply that the solution $\left(S,I,H\right)=\left(S(t),I(t),H(t)\right)$ remains in $\mathbb{R}_+^3$ for all time $t\geqslant 0$.

Next, letting $N=N(t)=S(t)+I(t)+H(t)$ for every $t\geqslant 0$ and adding the three equations in \eqref{eq:model3eqns}, one obtains
$$\frac{\dN}{\dt}=\lambda-\mu_1S-\mu_2I-\mu_2H\leqslant\lambda-\mu N,$$
where $\mu=\min\left\{\mu_1,\mu_2\right\}$. Multiplying both sides by $\e^{\mu t}$, one finds that the above inequality is equivalent to
$$\frac{\de}{\dt}\left[\e^{\mu t}N(t)\right]\leqslant \frac{\de}{\dt}\left[\frac{\lambda}{\mu}\e^{\mu t}-\frac{\lambda}{\mu}+N(0)\right].$$
Thus, at every time $t\geqslant 0$, the slope of the function $\e^{\mu t}N(t)$ is bounded above by that of the function $\left(\lambda/\mu\right)\e^{\mu t}-\lambda/\mu+N(0)$. Since both functions evaluate to $N(0)$ at $t=0$, this means that for every $t\geqslant 0$ we have
$$\e^{\mu t}N(t)\leqslant \frac{\lambda}{\mu}\e^{\mu t}-\frac{\lambda}{\mu}+N(0),$$
which means that
$$N(t)\leqslant \frac{\lambda}{\mu}+\left[N(0)-\frac{\lambda}{\mu}\right]\e^{-\mu t}.$$
Since the right-hand side approaches $\lambda/\mu$ as $t\to\infty$, we have the following theorem.\smallskip

\begin{theorem}
The sets $\mathbb{R}_+^3$ and
$$\left\{\left(S,I,H\right)\in\mathbb{R}_+^3:S+I+H\leqslant \frac{\lambda}{\min\left\{\mu_1,\mu_2\right\}}\right\}$$
are positively invariant under the model \eqref{eq:model3eqns}. Moreover, every solution of the model \eqref{eq:model3eqns}
associated to an initial condition in $\mathbb{R}_+^3$ is bounded.
\end{theorem}\smallskip

\subsection{Equilibria and basic reproduction number}\label{subsec:equilibria}

Let us now compute the equilibria of the model \eqref{eq:model3eqns} over its positively invariant domain $\mathbb{R}_+^3$. These are the solutions in $\mathbb{R}_+^3$ of the system
\begin{equation}\label{eq:FPsystem}
\left\{\begin{array}{rcl}
\displaystyle\lambda-\beta SI +\alpha_1 H +\alpha_2 I-\mu_1 S &\!\!\!=\!\!\!& \displaystyle0,\\
\displaystyle\beta SI -\left(\alpha_2+\gamma+\mu_2\right)I &\!\!\!=\!\!\!& \displaystyle0,\\
\displaystyle\gamma I-\left(\alpha_1+\mu_2\right)H &\!\!\!=\!\!\!& \displaystyle0.
\end{array}\right.
\end{equation}
The second equation can be written as
$$I\left(\beta S - \alpha_2-\gamma-\mu_2\right)=0,$$
which leads us to consider two cases, namely, the case $I=0$ and the case $S=\left(\alpha_2+\gamma+\mu_2\right)/\beta$.

\subsubsection{The disease-free equilibrium}

In the case of $I=0$, the system's third equation implies $H=0$, and the system's first equation subsequently implies $S=\lambda/\mu_1$. Consequently, $\mathbf{e}_{\DF}=\left(S_{\DF},I_{\DF},H_{\DF}\right)$, where
$$S_{\DF}=\frac{\lambda}{\mu_1},\qquad I_{\DF}=0,\qquad \text{and}\qquad H_{\DF}=0,$$
is the model's disease-free equilibrium. Notice that this equilibrium exists in the domain $\mathbb{R}_+^3$ for all parameter values.

\subsubsection{The basic reproduction number}

Before proceeding with the second case, let us compute the model's basic reproduction number. For this purpose, we employ the next-generation approach of van den Driessche and Watmough \cite{DriesscheWatmough,Martcheva}. The model's infected compartments are those of infected and hospitalised individuals, the numbers of which evolve via the equations
$$\frac{\dI}{\dt}=\mathcal{F}_I - \mathcal{V}_I,\qquad\text{where}\qquad \mathcal{F}_I:=\beta SI\quad\text{and}\quad \mathcal{V}_I=\left(\alpha_2+\gamma+\mu_2\right)I,$$
and
$$\frac{\dH}{\dt}=\mathcal{F}_H-\mathcal{V}_H,\qquad\text{where}\qquad \mathcal{F}_H=0\quad\text{and}\quad -\gamma I+\left(\alpha_1+\mu_2\right)H.$$
We thus define the matrices
$$\mathbf{F}:=\left(\begin{array}{cc}
\dfrac{\partial\mathcal{F}_I}{\partial I}\left(\mathbf{e}_{\DF}\right) & \dfrac{\partial\mathcal{F}_I}{\partial H}\left(\mathbf{e}_{\DF}\right)\\[0.3cm]
\dfrac{\partial\mathcal{F}_H}{\partial I}\left(\mathbf{e}_{\DF}\right) & \dfrac{\partial\mathcal{F}_H}{\partial H}\left(\mathbf{e}_{\DF}\right)
\end{array}\right)=\left(\begin{array}{cc}
\dfrac{\beta\lambda}{\mu_1} & 0\\[0.2cm]
0 & 0\end{array}\right)$$
and
$$\mathbf{V}:=\left(\begin{array}{cc}
\dfrac{\partial\mathcal{V}_I}{\partial I}\left(\mathbf{e}_{\DF}\right) & \dfrac{\partial\mathcal{V}_I}{\partial H}\left(\mathbf{e}_{\DF}\right)\\[0.3cm]
\dfrac{\partial\mathcal{V}_H}{\partial I}\left(\mathbf{e}_{\DF}\right) & \dfrac{\partial\mathcal{V}_H}{\partial H}\left(\mathbf{e}_{\DF}\right)
\end{array}\right)=\left(\begin{array}{cc}
\alpha_2+\gamma+\mu_2 & 0\\
-\gamma & \alpha_1+\mu_2
\end{array}\right),$$
so that the model's next-generation matrix reads
$$\mathbf{F}\mathbf{V}^{-1}=\left(\begin{array}{cc}
\dfrac{\beta\lambda}{\mu_1\left(\alpha_2+\gamma+\mu_2\right)} & 0\\[0.2cm]
0&0
\end{array}\right).$$
It follows that the model's basic reproduction number, which is the spectral radius of the above next-generation matrix, is given by
$$\cR_0=\frac{\beta\lambda}{\mu_1\left(\alpha_2+\gamma+\mu_2\right)}.$$
We are now ready to proceed with our second case.

\subsubsection{The endemic equilibrium}

In the second case, $S=\left(\alpha_2+\gamma+\mu_2\right)/\beta$, substituting this and the system's third equation into the system's first equation gives
$$\lambda-\left(\alpha_2+\gamma+\mu_2\right)I + \frac{\alpha_1\gamma}{\alpha_1+\mu_2}I+\alpha_2I-\frac{\mu_1\left(\alpha_2+\gamma+\mu_2\right)}{\beta}=0.$$
Multiplying both sides by $\cR_0/\lambda$ simplifies the equation to
$$\cR_0 - 1 = \frac{\mu_2\left(\alpha_1+\gamma+\mu_2\right)}{\lambda\left(\alpha_1+\mu_2\right)}\cR_0 I,$$
which gives
$$I=\frac{\lambda\left(\alpha_1+\mu_2\right)}{\mu_2\left(\alpha_1+\gamma+\mu_2\right)}\left(1-\frac{1}{\cR_0}\right).$$
The system's third equation then gives
$$H=\frac{\lambda\gamma}{\mu_2\left(\alpha_1+\gamma+\mu_2\right)}\left(1-\frac{1}{\cR_0}\right).$$
Therefore, $\mathbf{e}_{\E}=\left(S_{\E},I_{\E},H_{\E}\right)$, where
\begin{align*}
S_{\E}&=\frac{\alpha_2+\gamma+\mu_2}{\beta},\\
I_{\E}&=\frac{\lambda\left(\alpha_1+\mu_2\right)}{\mu_2\left(\alpha_1+\gamma+\mu_2\right)}\left(1-\frac{1}{\cR_0}\right),\\
H_{\E}&=\frac{\lambda\gamma}{\mu_2\left(\alpha_1+\gamma+\mu_2\right)}\left(1-\frac{1}{\cR_0}\right),
\end{align*}
is the model's endemic equilibrium. Clearly, this equilibrium exists in the domain $\mathbb{R}_+^3$ if and only if $\cR_0\geqslant 1$.\medskip

The following theorem summarises our results in this subsection.\smallskip

\begin{theorem}
The model \eqref{eq:model3eqns} over the positively invariant domain $\mathbb{R}_+^3$ has the basic reproduction number
\begin{equation}\label{eq:R0}
\cR_0=\frac{\beta\lambda}{\mu_1\left(\alpha_2+\gamma+\mu_2\right)}
\end{equation}
and two equilibria, namely, the disease-free equilibrium
\begin{equation}\label{eq:DFE}
\mathbf{e}_{\DF}=\left(\frac{\lambda}{\mu_1},\,0,\,0\right),
\end{equation}
which exists for all parameter values, and the endemic equilibrium
\begin{equation}\label{eq:EE}
\mathbf{e}_{\E}=\left(\frac{\alpha_2+\gamma+\mu_2}{\beta},\,\frac{\lambda\left(\alpha_1+\mu_2\right)}{\mu_2\left(\alpha_1+\gamma+\mu_2\right)}\left(1-\frac{1}{\cR_0}\right),\,\frac{\lambda\gamma}{\mu_2\left(\alpha_1+\gamma+\mu_2\right)}\left(1-\frac{1}{\cR_0}\right)\right),
\end{equation}
which exists if and only if $\cR_0\geqslant 1$.
\end{theorem}\smallskip

\subsection{Stability of equilibria}\label{subsec:stability}

Let us now relate the local asymptotic stability of the disease-free and endemic equilibria of the model \eqref{eq:model3eqns} to the model's basic reproduction number. For this purpose, we first compute the Jacobian $\mathbf{J}(S,I,H)$ of the model \eqref{eq:model3eqns}. Letting
$$\left\{\begin{array}{rcl}
f(S,I,H) &\!\!\!=\!\!\!& \displaystyle\lambda-\beta SI +\alpha_1 H +\alpha_2 I-\mu_1 S,\\
g(S,I,H) &\!\!\!=\!\!\!& \displaystyle\beta SI -\left(\alpha_2+\gamma+\mu_2\right)I,\\
h(S,I,H) &\!\!\!=\!\!\!& \displaystyle\gamma I-\left(\alpha_1+\mu_2\right)H,
\end{array}\right.$$
one obtains that
$$\mathbf{J}(S,I,H)=\left(\begin{array}{ccc}
\dfrac{\partial f}{\partial S} & \dfrac{\partial f}{\partial I} & \dfrac{\partial f}{\partial H}\\[0.3cm]
\dfrac{\partial g}{\partial S} & \dfrac{\partial g}{\partial I} & \dfrac{\partial g}{\partial H}\\[0.3cm]
\dfrac{\partial h}{\partial S} & \dfrac{\partial h}{\partial I} & \dfrac{\partial h}{\partial H}
\end{array}\right)=\left(\begin{array}{ccc}
-\beta I-\mu_1 & -\beta S + \alpha_2 & \alpha_1\\
\beta I & \beta S -\alpha_2-\gamma -\mu_2 & 0\\
0&\gamma&-\alpha_1-\mu_2
\end{array}\right).$$
The model's Jacobian near the disease-free equilibrium $\mathbf{e}_{\DF}$,
$$\mathbf{J}\left(\mathbf{e}_{\DF}\right)=\left(\begin{array}{ccc}
-\mu_1 & -\beta\lambda/\mu_1+\alpha_2 & \alpha_1\\
0 & \left(\alpha_2+\gamma+\mu_2\right)\left(\cR_0-1\right) & 0\\
0&\gamma&-\alpha_1-\mu_2
\end{array}\right),$$
has the characteristic polynomial
$$\left|r\mathbf{I}-\mathbf{J}\left(\mathbf{e}_{\DF}\right)\right|=\left(r+\mu_1\right)\left(r+\alpha_1+\mu_2\right)\left(r-\left(\alpha_2+\gamma+\mu_2\right)\left(\cR_0-1\right)\right).$$
The eigenvalues, therefore, are
$$r_1=-\mu_1<0,\qquad r_2=-\alpha_1-\mu_2<0,\qquad\text{and}\qquad r_3=\left(\alpha_2+\gamma+\mu_2\right)\left(\cR_0-1\right).$$
By \cite[Thm.\ 4.6]{Robinson}, we have proved the following theorem.\smallskip

\begin{theorem}
The disease-free equilibrium of the model \eqref{eq:model3eqns} is locally asymptotically stable if $\cR_0<1$, is unstable if $\cR_0>1$, and is non-hyperbolic if $\cR_0=1$.
\end{theorem}\smallskip

On the other hand, the model's Jacobian near the endemic equilibrium $\mathbf{e}_{\E}$ reads
$$\mathbf{J}\left(\mathbf{e}_{\E}\right)=\left(\begin{array}{ccc}
-\dfrac{\beta\lambda\left(\alpha_1+\mu_2\right)}{\mu_2\left(\alpha_1+\gamma+\mu_2\right)} \left(1-\dfrac{1}{\cR_0}\right)-\mu_1 & -\gamma-\mu_2 & \alpha_1\\[0.3cm]
\dfrac{\beta\lambda\left(\alpha_1+\mu_2\right)}{\mu_2\left(\alpha_1+\gamma+\mu_2\right)}\left(1-\dfrac{1}{\cR_0}\right) & 0 & 0\\[0.3cm]
0&\gamma & -\alpha_1-\mu_2
\end{array}\right).$$
Applying the substitution $\cR_0=\mathcal{S}+1$, one finds, through direct computation, that
$$\left|r\mathbf{I}-\mathbf{J}\left(\mathbf{e}_{\E}\right)\right|=r^3 + \mathcal{A}_1r^2 + \mathcal{A}_2 r + \mathcal{A}_3,$$
where
\begin{align*}
\cA_1 &= \frac{1}{\mu_2\left(\alpha_1+\gamma+\mu_2\right)\left(\cS+1\right)}\biggl[\left(\left(\alpha_1^2+\left(\gamma+\mu_1\right)\alpha_1+\beta\lambda+\mu_1\gamma\right)\cS+\left(\alpha_1+\mu_1\right)\left(\alpha_1+\gamma\right)\right)\mu_2\biggr.\\
&\phantom{=} \biggl.+\left(\cS+1\right)\mu_2^2\left(\mu_2+\gamma+2\alpha_1+\mu_1\right)+\alpha_1\beta\lambda\cS\biggr],\\
\cA_2 &= \frac{\left(\alpha_1+\mu_2\right)\left(\left(\cS+1\right)\mu_1\mu_2^2+\left(\left(\left(\alpha_1+\gamma\right)\mu_1+2\beta\lambda\right)\cS+\left(\alpha_1+\gamma\right)\mu_1\right)\mu_2+\cS\beta\lambda\left(\alpha_1+\gamma\right)\right)}{\mu_2\left(\alpha_1+\gamma+\mu_2\right)\left(\cS+1\right)},\\
\cA_3 &= \frac{\beta\lambda\left(\alpha_1+\mu_2\right)\cS}{\cS+1}.
\end{align*}
Suppose that $\cR_0>1$. Then, clearly, $\cS=\cR_0-1>0$, and so $\cA_1>0$ and $\cA_3>0$. Furthermore, direct computation shows that $\mu_2^2\left(\alpha_1+\gamma+\mu_2\right)^2\left(\cS+1\right)^2\left(\cA_1\cA_2-\cA_3\right)$ is a polynomial whose terms in the completely expanded form are all positive. Consequently, $\cA_1\cA_2-\cA_3>0$. By the Routh-Hurwitz criterion \cite[sec.\ 4.5]{Allen}, this implies that all the roots of $\left|r\mathbf{I}-\mathbf{J}\left(\mathbf{e}_{\E}\right)\right|$ have negative real parts. By \cite[Thm.\ 4.6]{Robinson}, we have proved the following theorem.\smallskip

\begin{theorem}
The endemic equilibrium of the model \eqref{eq:model3eqns} is locally asymptotically stable if $\cR_0>1$.
\end{theorem}\smallskip

\subsection{Design of insurance}\label{subsec:design}

Let us now design our health insurance based on the SIH-type model \eqref{eq:model3eqns}. To facilitate the design, we first complement the model with the two additional equations
\begin{equation}\label{eq:deaths}
\frac{\dD}{\dt}=\mu_1 S\qquad\text{and}\qquad 
\frac{\dDt}{\dt} =\mu_2 I+\mu_2H,
\end{equation}
which govern the time-evolution of the number $D=D(t)$ of natural deaths and the number $D^\ast=D^\ast(t)$ of deaths due to the disease any given time $t\geqslant 0$. Clearly, these numbers increase monotonically over time. Furthermore, recall that the unit for time in our model is chosen to be months. To begin our insurance design, let us specify a period for the insurance's availability, namely, a time-interval $[0,T]$, for some $T>0$. For every $t\in\{1,\ldots,T\}$, by month $t$ we mean the time-subinterval $[t-1,t]$.

\begin{table}\centering
	\renewcommand{\arraystretch}{1.9}\renewcommand{\tabcolsep}{6pt}
		\scalebox{0.925}{\begin{tabular}{|c|c|c|c|c|} 
			\hline
			Parameter & Description & Unit & 
\begin{minipage}{1.9cm}\centering\smallskip
Value for simulation\medskip
\end{minipage} & Source \\
			\hhline{|=|=|=|=|=|}
			$T$ & length of insurance availability period & month & 500 & simulated\\\hline
			$\Delta t$ & step size & month & 0.05000 & simulated\\\hline
			$i$ & monthly interest rate & -- & 0.00233 & \cite{Finance}\\\hline
			$\omega$ & \begin{minipage}{6.8cm}\centering\smallskip premium surcharge percentage allocated to operational costs\medskip\end{minipage} & -- & 0.10000 & simulated\\\hline
			$\varphi$ & \begin{minipage}{6.8cm}\centering\smallskip premium surcharge percentage allocated to profit\medskip\end{minipage} & -- & 0.05000 & simulated\\\hline
			$\cB_H$ &  monthly benefit for hospitalisation & \$/individual & 2,000 & simulated\\\hline
			$\cB_{D}$ & one-time benefit for natural death & \$/individual & 40,000 & simulated\\\hline
			$\cB_{D^\ast}$ & one-time benefit for death by disease & \$/individual & 50,000 & simulated\\\hline
			$D(0)$ & initial number of natural deaths & individual & 0 & simulated\\\hline
			$D^\ast(0)$ & initial number of deaths by disease & individual & 0 & simulated\\\hline
		\end{tabular}}\smallskip
		\caption{\label{tab:parameters2}Discretisation and insurance-related parameters used in our model and their values selected for our numerical simulations (section \ref{section:numerical}).}
\end{table}

Assume that for every $t\in\{1,\ldots,T\}$, the insurer is to receive from each individual who is either susceptible or infected at the beginning of the month a net premium of amount $\cPnet$. Assuming a constant monthly interest rate of $i$ and employing the standard notation $v=1/\left(1+i\right)$ for discounted cash flow, the present value of the total net premium over the stated period is thus given by
$$\Sigma\cPnet =\cPnet\sum_{t=0}^{T-1}v^t \left[S(t)+I(t)\right].$$
Next, assume that at the end of every month $t\in\{1,\ldots,T\}$, the insurer is to grant a benefit of amount $\cB_H$ to each individual who is currently hospitalised, a benefit of amount $\cB_D$ to each individual who dies naturally within the time-interval $(t-1,t]$, and a benefit of amount $\cB_{D^\ast}$ to each individual who dies due to the disease within the time-interval $(t-1,t]$. It follows that the present value of the total benefit granted over the stated period is given by
$$\Sigma\cB=\sum_{t=1}^T v^t \cB_H H(t) + \sum_{t=0}^{T-1} v^{t+1} \cB_D \left[D(t+1)-D(t)\right] + \sum_{t=0}^{T-1} v^{t+1} \cB_{D^\ast} \left[D^\ast(t+1)-D^\ast(t)\right].$$
The equivalence principle equation $\Sigma\cPnet=\Sigma\cB$ then gives the following expression for the insurance's monthly net premium:
$$\cPnet=\Sigma\cB\left[\sum_{t=0}^{T-1} v^t \left[S(t)+I(t)\right]\right]^{-1}.$$

Next, let us assume that the insurer adds surcharges of percentages $\omega$ and $\varphi$ to the monthly net premium $\cPnet$ to be allocated to operational costs and profit, respectively. The insurance's monthly gross premium to be paid by susceptible and infected individuals is thus given by $$\cPgross=\left(1+\omega+\varphi\right)\cPnet.$$ Consequently, the present value of the insurer's total profit at the end of every month $t\in\{1,\ldots,T\}$ is given by
$$\Pi(t) = \Sigma \cPgross(t) -\Sigma \cC(t) -\Sigma\cB(t),$$
where
$$\Sigma \cPgross(t) = \left(1+\omega+\varphi\right)\cPnet\sum_{\tau=0}^{t-1}v^\tau \left[S(\tau)+I(\tau)\right]$$
is the present value of the total gross premium the insurer receives until the beginning of month $t$,
$$\Sigma \cC(t) = \omega\,\cPnet\sum_{\tau=0}^{t-1}v^\tau \left[S(\tau)+I(\tau)\right]$$
is the present value of the total operational costs the insurer incurs until the beginning of month $t$, and
$$\Sigma\cB(t)= \sum_{\tau=1}^t v^\tau \cB_H H(\tau) +\sum_{\tau=0}^{t-1} v^{\tau+1} \cB_D \left[D(\tau+1)-D(\tau)\right] + \sum_{\tau=0}^{t-1} v^{\tau+1} \cB_{D^\ast} \left[D^\ast(\tau+1)-D^\ast(\tau)\right]$$
is the present value of the total benefit the insurer grants until the end of month $t$. For completeness, let us define $\Pi(0)=0$.

Finally, let us define
$$\Pi_{\min} = \min_{t\in\{0,\ldots,T\}}\Pi(t)\qquad\text{and}\qquad t_{\min} = \underset{t\in\{0,\ldots,T\}}{\argmin}\Pi(t), $$
assuming uniqueness for the latter. In the case $\Pi_{\min}\geqslant 0$, the insurer never suffers a loss throughout the stated period. Let us now consider the case $\Pi_{\min}<0$, where the insurer suffers a loss of maximum amount $-\Pi_{\min}=-\Pi\left(t_{\min}\right)$ over the stated period. In order to prevent this loss, it is necessary for the insurer to allocate a start-up capital of minimum amount
$$\Gamma=-\Pi_{\min} v^{t_{\min}}$$
in the beginning of the first month. Assuming the allocation of such a capital, the present value of the insurer's financial asset at the end of every month $t\in\{1,\ldots,T\}$,
$$\Gamma(t)=\Gamma+\Pi(t),$$ 
satisfies $\Gamma(t)\geqslant 0$ with equality if and only if $t=t_{\min}$. Since the present value of the insurer's total profit at the end of the insurance's availability period, i.e., at the end of month $T$, is precisely $\Pi(T)$, we see that the percentage of this total profit with respect to the start-up capital is given by
$$\pi=\frac{\Pi(T)}{\Gamma}\times 100\%.$$
A summary of the parameters involved in our insurance design is provided by Table \ref{tab:parameters2}.

\section{Discretisation}\label{sec:discretisation}

To prepare a setting for our numerical simulation, let us now construct a discrete version of our continuous model \eqref{eq:model3eqns} using the forward Euler method (subsection \ref{subsec:discrete}). It shall be immediately apparent that the equilibria of our discrete model are precisely those of our continuous model. We establish sufficient conditions for the equilibria's stability, involving not only the model's basic reproduction number but also the discretisation step size (subsection \ref{subsec:stabilitydiscrete}). Finally, we develop discrete analogues of the previously-formulated financial quantities associated to our health insurance (subsection \ref{subsec:insurancediscrete}).

\subsection{Discrete model and equilibria}\label{subsec:discrete}

Let us discretise the model \eqref{eq:model3eqns} using the forward Euler method \cite[sec.\ 22.3]{KongSiauwBayen}, with an arbitrary step size $\Delta t>0$ (Table \ref{tab:parameters2}). Accordingly, for every non-negative integer $n$, we introduce the time-step
$$t_n=n\Delta t$$
and the time-evolving variables
$$\tS_n\approx S\left(t_n\right),\qquad \tI_n\approx I\left(t_n\right),\qquad\text{and}\qquad \tH_n\approx H\left(t_n\right)$$
governed by the recursion
\begin{equation}\label{eq:modeldiscrete}
\left\{\begin{array}{rcl}
\tS_{n+1} &\!\!\!=\!\!\!& \tS_n +  \left(\lambda-\beta \tS_n\tI_n +\alpha_1 \tH_n +\alpha_2 \tI_n-\mu_1 \tS_n\right)\Delta t,\\
\tI_{n+1} &\!\!\!=\!\!\!& \tI_n +\left(\beta \tS_n\tI_n -\left(\alpha_2+\gamma+\mu_2\right)\tI_n\right) \Delta t ,\\
\tH_{n+1} &\!\!\!=\!\!\!& \tH_n+\left(\gamma \tI_n-\left(\alpha_1+\mu_2\right)\tH_n\right)\Delta t,
\end{array}\right.
\end{equation}
along with the initial condition $\left(\tS_0,\tI_0,\tH_0\right)=\left(S(0),I(0),H(0)\right)$.

Clearly, the equilibria in $\mathbb{R}_+^3$ of the discrete model \eqref{eq:modeldiscrete} are the solutions in $\mathbb{R}_+^3$ of the system \eqref{eq:FPsystem}, namely, the equilibria $\mathbf{e}_{\DF}$ and $\mathbf{e}_{\E}$ of the continuous model \eqref{eq:model3eqns}, given by \eqref{eq:DFE} and \eqref{eq:EE}. In the upcoming subsection, we study the stability of these equilibria.

\subsection{Stability of discrete model's equilibria}\label{subsec:stabilitydiscrete}

Letting
$$\left\{\begin{array}{rcl}
\overline{f}(\tS,\tI,\tH) &\!\!\!=\!\!\!& \displaystyle \tS+\left(\lambda-\beta \tS\,\tI +\alpha_1 \tH +\alpha_2 \tI-\mu_1 \tS\right)\Delta t,\\
\overline{g}(\tS,\tI,\tH) &\!\!\!=\!\!\!& \displaystyle \tI + \left(\beta \tS\,\tI -\left(\alpha_2+\gamma+\mu_2\right)\tI\right)\Delta t,\\
\overline{h}(\tS,\tI,\tH) &\!\!\!=\!\!\!& \displaystyle \tH + \left(\gamma \tI-\left(\alpha_1+\mu_2\right)\tH\right)\Delta t,
\end{array}\right.$$
the Jacobian of the discrete model \eqref{eq:modeldiscrete} is given by
\begin{align*}
\overline{\mathbf{J}}\left(\tS,\tI,\tH\right)&=\left(\begin{array}{ccc}
\dfrac{\partial \overline f}{\partial \tS} & \dfrac{\partial \overline f}{\partial \tI} & \dfrac{\partial \overline f}{\partial \tH}\\[0.3cm]
\dfrac{\partial \overline g}{\partial \tS} & \dfrac{\partial \overline g}{\partial \tI} & \dfrac{\partial \overline g}{\partial \tH}\\[0.3cm]
\dfrac{\partial \overline h}{\partial \tS} & \dfrac{\partial \overline h}{\partial \tI} & \dfrac{\partial \overline h}{\partial \tH}
\end{array}\right)\\
&=\left(\begin{array}{ccc}
1-\left(\beta \tI+\mu_1\right)\Delta t & \left(-\beta \tS + \alpha_2\right)\Delta t &\alpha_1 \Delta t\\
 \beta \tI\Delta t & 1-\left(\alpha_2+\gamma +\mu_2-\beta \tS\right)\Delta t & 0\\
0&\gamma\Delta t&1-\left(\alpha_1+\mu_2\right)\Delta t
\end{array}\right).
\end{align*}
Evaluating this near the disease-free equilibrium $\mathbf{e}_{\DF}$ gives the matrix
$$\overline{\mathbf{J}}\left(\mathbf{e}_{\DF}\right)=\left(\begin{array}{ccc}
1- \mu_1\Delta t & \left(-\beta\lambda/\mu_1 + \alpha_2\right)\Delta t & \alpha_1\Delta t\\
0 & 1-\left(\alpha_2+\gamma+\mu_2\right)\left(1-\cR_0\right)\Delta t & 0\\
0&\gamma\Delta t&1-\left(\alpha_1+\mu_2\right)\Delta t
\end{array}\right),$$
%Characteristic polynomial
%$$\left|r\mathbf{I}-\overline{\mathbf{J}}\left(\mathbf{e}_{\DF}\right)\right|=\left(r-\left(1-\mu_1\Delta t\right)\right)\left(r-\left(1-\left(\alpha_2+\gamma+\mu_2\right)\left(1-\cR_0\right)\Delta t\right)\right)\left(r-\left(1-\left(\alpha_1+\mu_2\right)\Delta t\right)\right).$$
possessing the eigenvalues
$$\overline r_1= 1-\mu_1\Delta t,\quad \overline r_2 = 1-\left(\alpha_2+\gamma+\mu_2\right)\left(1-\cR_0\right)\Delta t,\quad \text{and}\quad \overline r_3=1-\left(\alpha_1+\mu_2\right)\Delta t.$$
By \cite[Theorem 12.3]{Robinson}, this leads to the following theorem.\smallskip

\begin{theorem}
The disease-free equilibrium of the discrete model \eqref{eq:modeldiscrete} is locally asymptotically stable if 
$$\cR_0<1\,\,\,\text{and}\,\,\, \Delta t<\min\left\{\frac{2}{\mu_1},\frac{2}{\left(\alpha_2+\gamma+\mu_2\right)\left(1-\cR_0\right)},\frac{2}{\alpha_2+\mu_2}\right\};$$
is unstable if
$$\cR_0>1\qquad\text{or}\qquad \cR_0<1\,\,\,\text{and}\,\,\, \Delta t>\min\left\{\frac{2}{\mu_1},\frac{2}{\left(\alpha_2+\gamma+\mu_2\right)\left(1-\cR_0\right)},\frac{2}{\alpha_2+\mu_2}\right\};$$
and is non-hyperbolic if
$$\cR_0=1\qquad\text{or}\qquad \Delta t\in\left\{\frac{2}{\mu_1},\frac{2}{\left(\alpha_2+\gamma+\mu_2\right)\left(1-\cR_0\right)},\frac{2}{\alpha_2+\mu_2}\right\}.$$
\end{theorem}\smallskip

Next, the discrete model's Jacobian near the endemic equilibrium $\mathbf{e}_{\E}$,
$$\overline{\mathbf{J}}\left(\mathbf{e}_{\E}\right) =\left(\begin{array}{ccc}
1-\left(\beta I_{\E}+\mu_1\right)\Delta t & \left(-\gamma-\mu_2\right)\Delta t & \alpha_1\Delta t \\
\beta I_{\E} \Delta t  & 1 & 0\\
0&\gamma\Delta t  & 1-\left(\alpha_1-\mu_2\right)\Delta t 
\end{array}\right),$$
has the characteristic polynomial
$$\bigl|r\mathbf{I}-\overline{\mathbf{J}}\left(\mathbf{e}_{\E}\right)\bigr|=r^3 + \overline{\cA}_1 r^2 + \overline{\cA}_2 r + \overline{\cA}_3,$$
where
\begin{align*}
\overline{\cA}_1 &= \left(\beta I_{\E}+\alpha_1+\mu_1+\mu_2\right)\Delta t-3,\\
\overline{\cA}_2 &= \left[\beta\left(\alpha_1+\gamma+2\mu_2\right)I_{\E} +\mu_1\left(\alpha_1+\mu_2\right)\right]\left(\Delta t\right)^2-2\left(\beta I_{\E}+\alpha_1+\mu_1+\mu_2\right)\Delta t+3,\\
\overline{\cA}_3 &= \beta I_{\E}\mu_2\left(\alpha_1+\gamma+\mu_2\right)\left(\Delta t\right)^3- \left[\beta I_{\E}\left(\alpha_1+\gamma+2\mu_2\right)+\mu_1\left(\alpha_1+\mu_2\right)\right]\left(\Delta t\right)^2 \\
&\quad + \left(\beta I_{\E}+\alpha_1+\mu_1+\mu_2\right)\Delta t-1.
\end{align*}
In the case of $\cR_0>1$, we have that $I_{\E}>0$, and so
$$1+\overline{\cA}_1+\overline{\cA}_2+\overline{\cA}_3= \beta I_{\E}\mu_2\left(\alpha_1+\gamma+\mu_2\right)\left(\Delta t\right)^3>0.$$
Consequently, the Schur-Cohn criterion \cite[sec.\ 5.1]{Elaydi} implies the following theorem.\smallskip

\begin{theorem}
Suppose that $\cR_0>1$. The endemic equilibrium of the discrete model \eqref{eq:modeldiscrete} is locally asymptotically stable if 
$$1-\overline{\cA}_1+\overline{\cA}_2-\overline{\cA}_3>0\qquad\text{and}\qquad \left|\overline{\cA}_2-\overline{\cA}_1\overline{\cA}_3\right|<1-\overline{\cA}_3^2.$$
\end{theorem}\smallskip

\subsection{Discrete analogues of insurance-related quantities}\label{subsec:insurancediscrete}

In this subsection, we develop the discrete analogues of the insurance-related quantities formulated in subsection \ref{subsec:design}. First, the forward Euler method discretises the equations \eqref{eq:deaths} to
$$\tD_{n+1}=\tD_n+\mu_1 \tS_n\Delta t\qquad\text{and}\qquad \tD^\ast_{n+1}=\tD^\ast_n +\left(\mu_2 \tI_n+\mu_2\tH_n\right)\Delta t,$$
respectively, where
$$\tD_n\approx D\left(t_n\right)\qquad \text{and}\qquad \tD^\ast_n \approx D^\ast\left(t_n\right),$$
with the initial condition $\bigr(\tD_0,\tD^\ast_0\bigr)=\left(D(0),D^\ast (0)\right)$.

For computational convenience, let us assume that the step size $\Delta t$ is chosen such that $1/\Delta t$ is a natural number. Letting $\overline{\cP}^{\text{net}}$ be an estimate for the insurance's monthly net premium $\cPnet$, one obtains the following estimate for the present value of the total net premium $\Sigma\cPnet$ over the period $[0,T]$:
$$\overline{\Sigma\mathcal{P}}^{\text{net}} =\overline{\cP} \sum_{t=0}^{T-1}v^t \left(\tS_{t/\Delta t}+\tI_{t/\Delta t}\right).$$
On the other hand, an estimate for the present value of the total benefit granted over the period is given by
$$\overline{\Sigma\cB}=\sum_{t=1}^T v^t \cB_H \tH_{t/\Delta t} + \sum_{t=0}^{T-1} v^{t+1} \cB_D \left[\tD_{(t+1)/\Delta t}-\tD_{t/\Delta t}\right]+ \sum_{t=0}^{T-1} v^{t+1} \cB_{D^\ast} \left[\tD^\ast_{(t+1)/\Delta t}-\tD^\ast_{t/\Delta t}\right].$$
The equivalence principle equation $\overline{\Sigma\mathcal{P}}^{\text{net}}=\overline{\Sigma\cB}$ thus gives the following expression for the insurance's monthly net premium:
$$\overline{\cP}^{\text{net}}=\overline{\Sigma\mathcal{B}}\left[\sum_{t=0}^{T-1}v^t \left(\tS_{t/\Delta t}+\tI_{t/\Delta t}\right)\right]^{-1}.$$

Consequently, $$\overline{\cP}^{\text{gross}}=\left(1+\omega+\varphi\right)\overline{\cP}^{\text{net}}$$ estimates the insurance's monthly gross premium $\cPgross$. Furthermore, for every $t\in\{1,\ldots,T\}$, estimates for the present value $\Sigma \cPgross(t)$ of the total gross premium the insurer receives until the beginning of month $t$, for the present value $\Sigma \cC(t)$ of the total operational costs the insurer incurs until the beginning of month $t$, and for the present value $\Sigma\cB(t)$ of the total benefit the insurer grants until the end of month $t$ are given by
\begin{align*}
\overline{\Sigma\cP}^{\text{gross}}_{t/\Delta t}&=\left(1+\omega+\varphi\right)\overline{\cP}^{\text{net}}\sum_{\tau=0}^{t-1}v^\tau \left(\tS_{\tau/\Delta t}+\tI_{\tau/\Delta t}\right),\\
\overline{\Sigma\cC}_{t/\Delta t}&=\omega\,\overline{\cP}^{\text{net}}\sum_{\tau=0}^{t-1}v^\tau \left(\tS_{\tau/\Delta t}+\tI_{\tau/\Delta t}\right),\\
\overline{\Sigma\cB}_{t/\Delta t}&=\sum_{\tau=1}^t v^\tau \cB_H \tH_{\tau/\Delta t} + \sum_{\tau=0}^{t-1} v^{\tau+1} \cB_D \left[\tD_{(\tau+1)/\Delta t}-\tD_{\tau/\Delta t}\right]+ \sum_{\tau=0}^{t-1} v^{\tau+1} \cB_{D^\ast}  \left[\tD^\ast_{(\tau+1)/\Delta t}-\tD^\ast_{\tau/\Delta t}\right],
\end{align*}
respectively. It follows that the present value $\Pi(t)$ of the insurer's total profit at the end of every month $t\in\{1,\ldots,T\}$ is estimated by
$$\overline{\Pi}_{t/\Delta t}=\overline{\Sigma\cP}^{\text{gross}}_{t/\Delta t} -\overline{\Sigma\cC}_{t/\Delta t} -\overline{\Sigma\cB}_{t/\Delta t}.$$
As before, for completeness, we define $\overline{\Pi}_0=0$. 

Next, the quantities
\begin{equation}\label{eq:minTP}
\overline{\Pi}_{\min} = \min_{t\in\{0,\ldots,T\}}\overline{\Pi}_{t/\Delta t}\qquad\text{and}\qquad \overline{t}_{\min} = \underset{t\in\{0,\ldots,T\}}{\argmin}\overline{\Pi}_{t/\Delta t}, 
\end{equation}
the latter defined assuming uniqueness, estimate $\Pi_{\min}$ and $t_{\min}$, respectively. If $\overline{\Pi}_{\min}\geqslant 0$, then no loss is experienced by the insurer throughout the insurance's availability period. Otherwise, the insurer experiences a loss of maximum amount $-\overline{\Pi}_{\min}=-\overline{\Pi}_{\overline{t}_{\min}/\Delta t}$, so that the minimum loss-preventing start-up capital to be allocated at the beginning of the insurance's availability period is given by
\begin{equation}\label{eq:capital}
\overline{\Gamma}=-\overline{\Pi}_{\min} v^{\overline{t}_{\min}},
\end{equation}
which estimates $\Gamma$. In the case of such a capital being allocated, the present value of the insurer's financial asset $\Gamma_t$ at the end of every month $t\in\{0,\ldots,T\}$ is estimated by
$$\overline{\Gamma}_{t/\Delta t}=\overline{\Gamma}+\overline{\Pi}_{t/\Delta t}.$$
Finally, since $\overline{\Pi}_{T/\Delta t}$ estimates the present value of the insurer's total profit $\Pi(T)$ at the end of month $T$, then 
\begin{equation}\label{eq:percent}
\overline{\pi}=\frac{\overline{\Pi}_{T/\Delta t}}{\overline{\Gamma}}\times 100\%
\end{equation}
estimates the percentage $\pi$ of the insurer's total profit with respect to the allocated start-up capital.

\section{Numerical simulation and sensitivity analysis}\label{section:numerical}

In this section, we employ our discrete model to simulate the spread of a generic infectious disease characterised by the two sets of parameter values provided by Tables \ref{tab:parameters1} and \ref{tab:parameters2}, which differ only in the value of the incidence coefficient: $\beta=0.001$ representing the subcritical case $\cR_0<1$, and $\beta=0.003$ representing the supercritical case $\cR_0>1$. These values are selected from a variety of sources, in the way we shall explain in subsection \ref{subsec:parameter}. In subsection \ref{subsec:numerical}, we shall visualise the time-evolution of the number of susceptible, infected, and hospitalised individuals as governed by our discrete model, and analyse the financial quantities associated to our health insurance, in each of the two cases. Finally, in subsection \ref{subsec:sensitivity}, we assess the sensitivity of four key quantities: the model's basic reproduction number $\cR_0$, the insurance's monthly gross premium $\overline{\cP}^{\text{gross}}$, the insurer's start-up capital $\overline{\Gamma}$, and the insurer's end-of-period total profit $\overline{\Pi}_{T/\Delta t}$, with respect to each of the model's parameters, again in each of the two cases.

\subsection{Selection of parameter values}\label{subsec:parameter}

Let us begin by justifying our selection of parameter values presented in Tables \ref{tab:parameters1} and \ref{tab:parameters2}.
\begin{enumerate}[leftmargin=1cm,itemsep=3pt]
\item[(i)] The initial numbers $S(0)$, $I(0)$, $H(0)$, $D(0)$, and $D^*(0)$ of susceptible individuals, infected individuals, hospitalised individuals, natural deaths, and deaths by disease are all chosen to simulate a population initially comprising $N(0)=3000$ individuals, with $S(0)=2999$ being susceptible and $I(0)=1$ being infected.

\item[(ii)] Our population's birth rate $\lambda$ is estimated as follows. First, since our population initially consists of $N(0)=3000$ individuals, and the global population is estimated to consist of 102 males for every 100 females \cite[p.\ 2]{Affairs}, we infer that there are approximately $100/\left(100+102\right)\cdot 3000\approx 1485$ females in our initial population. Next, we note that the global number of births per female and the global life expectancy are recently estimated to be 2.5 individuals per female \cite{UN} and 73.4 years \cite{WHO}, respectively. Therefore, we make the estimation $\lambda\approx \left(1485 \cdot 2.5\right)/\left(73.4\cdot 12\right) \approx 4.21492$ individuals per month.

\item[(iii)] For the treatment coefficient $\alpha_1$ of hospitalised infected individuals, we refer to the work of Bettger et al.\ on recovery patterns of such individuals \cite{BettgerEtAl}, which reported that 60\% of such individuals show significant improvement within a single year of treatment, leading to the estimation $\alpha_1\approx 0.6/12=0.05000$ per month. For the treatment coefficient $\alpha_2$ of non-hospitalised infected individuals, we assume that the substandard treatment received by such individuals mitigates their typically milder symptoms to an extent which justifies assigning them the same treatment coefficient, namely, $\alpha_2\approx 0.05000$ per month.

\item[(iv)] Our estimate for the hospitalisation coefficient $\gamma$ of infected individuals originate from the work of Nichol et al.\ \cite{NicholEtAl}, which observed that during non-influenza seasons, unvaccinated pneumonia and influenza patients are hospitalised at the rate of 55 per 1000 person-years. This leads to the estimation $\gamma\approx \left(55/100\right)\cdot 12=0.66000$ per month.

\item[(v)] The 2015 static mortality table of the Internal Revenue Service \cite{SOA} provides the estimated probabilities that a person of age from 1 to 119 in the United States dies within a year, which average to approximately 0.08942. Accordingly, we assume that $\mu_1\approx 0.08942/12\approx 0.00745$ of the individuals in our population die naturally every month.

\item[(vi)] During the 7.5-month period from April 1 to November 15, 2020, Malhotra et al.\ \cite{MalhotraEtAl} recorded in a hospital in New Delhi that approximately 13.72\% of laboratory-confirmed patients with COVID-19 died. Based on this data, we assume that $\mu_2\approx 13.72\%/7.5\approx 0.01829$ of the individuals in our population die due to the disease every month.

\item[(vii)] Our constant monthly interest rate $i$ is estimated from the average adjusted close percentages of the 30-year treasury yield rate, which serves as as a benchmark for risk-free rates \cite{Finance}. This average, $2.79805\%$, translates to a monthly rate of $i\approx 2.79805\%/12\approx 0.00233$.

\item[(viii)] The values of all other parameters are simulated. In particular, for the incidence coefficient $\beta$, we select two different values to represent a disease-free and an endemic scenario: 0.00100 and 0.00300, that is, 1 and 3 per 1000 person-months, respectively.
\end{enumerate}

\subsection{Numerical solution and financial analysis}\label{subsec:numerical}

In each of the two scenarios ---a disease-free scenario and an endemic scenario--- described by our two sets of parameter values provided in Tables \ref{tab:parameters1} and \ref{tab:parameters2}, let us visualise the solution of our discrete model \eqref{eq:modeldiscrete} and calculate the financial quantities associated to our insurance.

\subsubsection{A disease-free scenario}

\begin{figure}
\includegraphics[height=5cm]{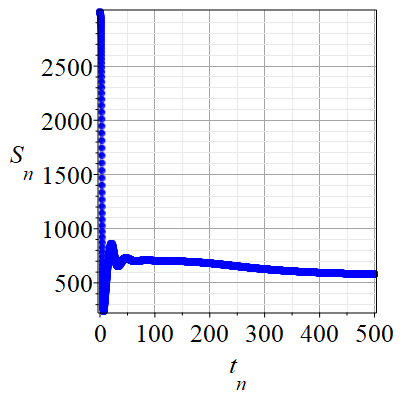}\,\,\includegraphics[height=5cm]{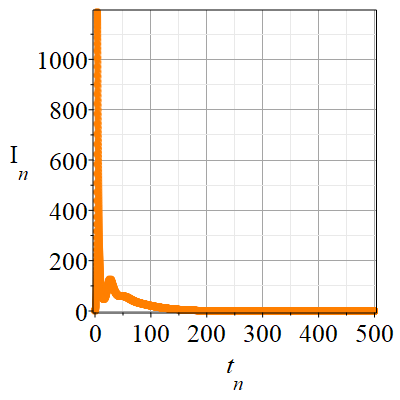}\,\,\includegraphics[height=5cm]{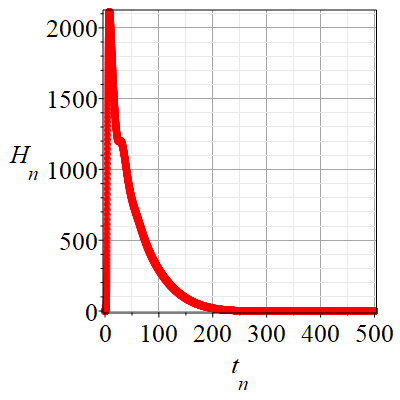} 
\caption{\label{fig:TEdf}The time-evolution of the number of susceptible (left), infected (middle), and hospitalised individuals (right), where the parameter values are given in Table \ref{tab:parameters1} with $\beta=0.00100$. In this case, $\cR_0\approx 0.77683<1$.}
\end{figure}

\begin{figure}
\includegraphics[height=5cm]{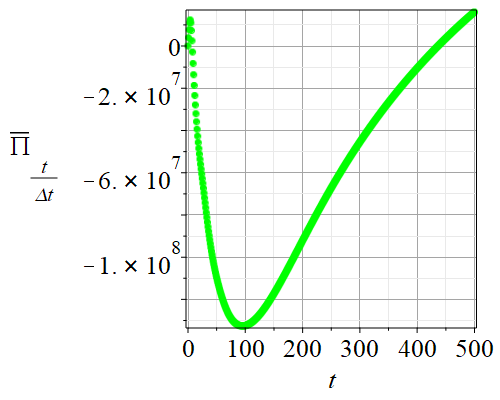}\,\,\,\includegraphics[height=5cm]{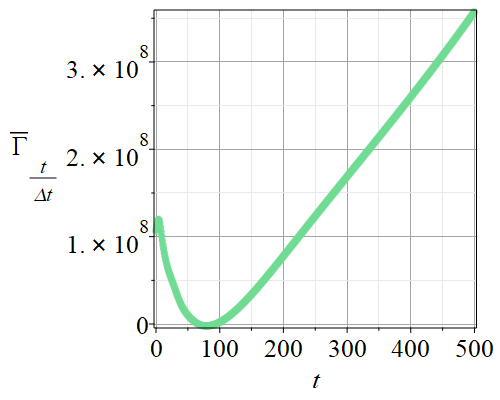} 
\caption{\label{fig:TEPdf}The time-evolution of the insurer's total profit (left) and of the insurer's total asset (right), where the parameter values are given in Table \ref{tab:parameters2} with $\beta=0.00100$. In this case, $\cR_0\approx 0.77683<1$.}
\end{figure}

We consider the set of parameter values provided by Tables \ref{tab:parameters1} and \ref{tab:parameters2}, with $\beta=0.00100$. For these parameter values, the basic reproduction number \eqref{eq:R0} evaluates to $\cR_0\approx 0.77683<1$. The time-evolution of the numbers of susceptible, infected, hospitalised individuals are visualised in Figure \ref{fig:TEdf}. Notice the convergence of these numbers to the coordinates $S_\infty$, $I_\infty$, $H_\infty$ of the disease-free equilibrium $\mathbf{D}_{\text{DF}}\approx\left(566, 0, 0\right)$ given by \eqref{eq:DFE}. Our insurance's monthly gross premium evaluates to $\overline{\cP}^{\text{gross}}\approx 1,738$ dollars, and the time-evolution of the present value of the insurer's total profit is visualised in Figure \ref{fig:TEPdf} (left). The minimum of these present values, $\overline{\Pi}_{\min}\approx -132,583,472$, is negative, and is achieved at $t=\overline{t}_{\min}=95$; see equation \eqref{eq:minTP}. This leads to the necessity of allocating a start-up capital of $\overline{\Gamma}\approx 106,284,546$ dollars, computed using equation \eqref{eq:capital}. Assuming the allocation of such a capital, the time-evolution of the present value of the insurer's total asset is visualised in Figure \ref{fig:TEPdf} (right). The present value of the insurer's total profit at the end of the insurance's availability period is estimated to be $\overline{\Pi}_{T/\Delta t}\approx 16,106,242$ dollars, which is $\overline{\pi}\approx 15.15389\%$ of the allocated start-up capital, by equation \eqref{eq:percent}. We summarise our results in this scenario in the third column of Table \ref{tab:numerical}.

\subsubsection{An endemic scenario}

\begin{figure}
\includegraphics[height=5cm]{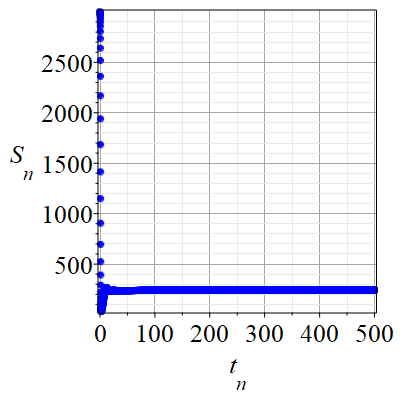}\,\,\includegraphics[height=5cm]{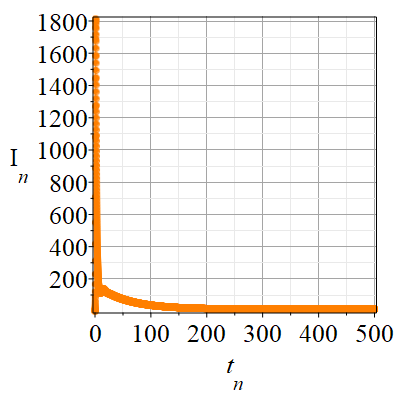}\,\,\includegraphics[height=5cm]{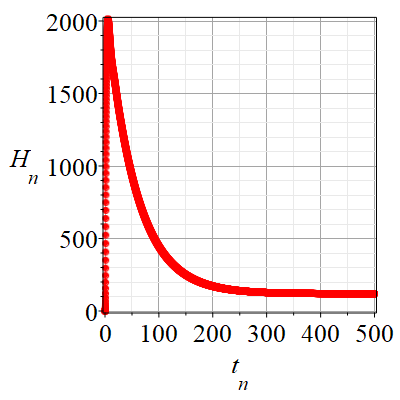} 
\caption{\label{fig:TEe}The time-evolution of the number of susceptible (left), infected (middle), and hospitalised individuals (right), where the parameter values are given in Table \ref{tab:parameters1} with $\beta=0.00300$. In this scenario, $\cR_0\approx 2.33050>1$.}
\end{figure}

\begin{figure}
\includegraphics[height=5cm]{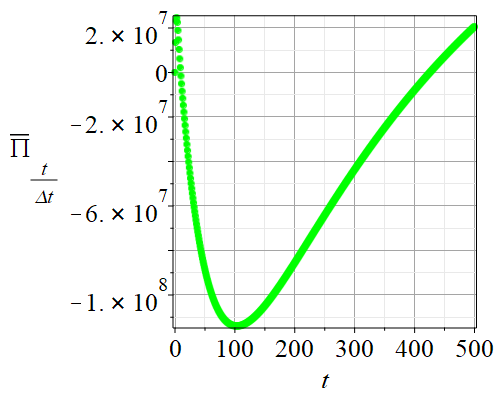}\,\,\,\includegraphics[height=5cm]{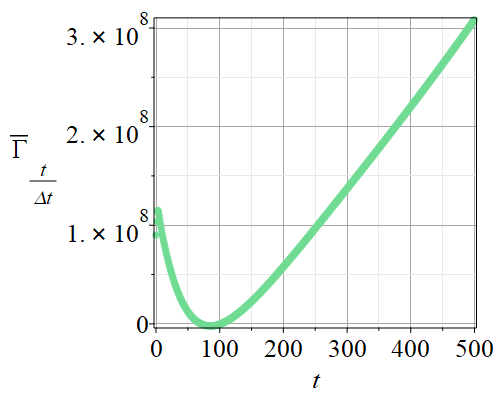} 
\caption{\label{fig:TEPe}The time-evolution of the insurer's total profit (left) and of the insurer's total asset (right), where the parameter values are given in Table \ref{tab:parameters2} with $\beta=0.00300$. In this scenario, $\cR_0\approx 2.33050>1$.}
\end{figure}

For a comparison, let us now present the numerical results obtained in our endemic scenario, characterised by the parameter values provided by Tables \ref{tab:parameters1} and \ref{tab:parameters2}, with $\beta=0.00300$. In this scenario, the model's basic reproduction number \eqref{eq:R0} takes the value of $\cR_0\approx2.33050>1$, and the time-evolution of the numbers of susceptible, infected, hospitalised individuals are visualised in Figure \ref{fig:TEe}. These numbers converge to the coordinates $S_\infty$, $I_\infty$, $H_\infty$ of the endemic equilibrium $\mathbf{D}_{\text{E}}\approx\left(243, 12, 119\right)$ prescribed by \eqref{eq:EE}. As expected, our insurance's monthly gross premium in this scenario is significantly higher, namely, $\overline{\cP}^{\text{gross}}\approx 5,338$ dollars. However, the present value of the insurer's total profit, whose time-evolution is visualised in Figure \ref{fig:TEPe} (left), achieves a slightly less negative minimum value of $\overline{\Pi}_{\min}\approx -113,944,943$, at $t=\overline{t}_{\min}=103$, leading to a lower minimum necessary start-up capital of $\overline{\Gamma}\approx 89,658,189$ dollars, by equations \eqref{eq:minTP} and \eqref{eq:capital}. With such a capital allocated, Figure \ref{fig:TEPe} (right) visualises the time-evolution of the present value of the insurer's total asset. Finally, the present value of the insurer's end-of-period total profit evaluates to $\overline{\Pi}_{T/\Delta t}\approx 20,590,132$ dollars, which is $\overline{\pi}\approx 22.96514\%$ of the allocated start-up capital, by equation \eqref{eq:percent}. Our results in this scenario are summarised in the fourth column of Table \ref{tab:numerical}. We thus see that, in the endemic scenario, despite the significantly higher value of monthly gross premium, the minimum start-up capital the insurer must allocate to avoid loss is lower, and this leads to a higher total profit percentage.

\begin{table}\renewcommand{\arraystretch}{1.9}\renewcommand{\tabcolsep}{4pt}
\scalebox{0.925}{\begin{tabular}{|c|c|c|c|c|}\hline
\multirow{2}{*}{Quantity} & \multirow{2}{*}{Description} & \multicolumn{2}{c|}{Value}& \multirow{2}{*}{Unit} \\\cline{3-4}
&& \begin{minipage}{2.25cm}\centering Disease-free\end{minipage} & \begin{minipage}{2.25cm}\centering Endemic\end{minipage}&\\\hhline{|=|=|=|=|=|}
$\beta$ & incidence coefficient & $0.00100$ & $0.00300$ & 1/(individual\,$\times$\,month)\\\hline
$\cR_0$ & basic reproduction number & $0.77683$ & $2.33050$ & --\\\hline
$S_\infty$ & \begin{minipage}{6cm}\centering\smallskip limiting number of susceptible individuals\medskip\end{minipage} & $566$ & $243$ & individual\\\hline
$I_\infty$ & \begin{minipage}{6cm}\centering\smallskip limiting number of non-hospitalised infected individuals\medskip\end{minipage} & $0$ & $12$ & individual\\\hline
$H_\infty$ & \begin{minipage}{6cm}\centering\smallskip limiting number of hospitalised individuals\medskip\end{minipage} & $0$ & $119$ & individual\\\hline
$\overline{\cP}^{\text{gross}}$ & monthly gross premium & $1,738$ & $5,338$ & \$/individual\\\hline
$\overline{\Pi}_{\min}$ & \begin{minipage}{6cm}\centering\smallskip minimum present value of total profit\medskip\end{minipage} & $-132,583,472$ & $-113,944,943$ & \$\\\hline
$\overline{t}_{\min}$ & \begin{minipage}{6cm}\centering\smallskip time at which minimum present value of total profit is achieved\medskip\end{minipage} & $95$ & $103$ & month\\\hline
$\overline{\Gamma}$ & \begin{minipage}{6cm}\centering\smallskip minimum loss-preventing start-up capital\medskip\end{minipage} & $106,284,546$ & $89,658,189$ & \$\\\hline
$\overline{\Pi}_{T/\Delta t}$ & \begin{minipage}{6cm}\centering\smallskip present value of end-of-period total profit\medskip\end{minipage} & $16,106,242$ & $20,590,132$ & \$\\\hline
$\overline{\pi}$ & total profit percentage& $15.15389\%$ & $22.96514\%$ & --\\\hline
\end{tabular}}\smallskip
\caption{\label{tab:numerical}Summary of different quantities in our two numerical scenarios.}
\end{table}

\subsection{Sensitivity analysis}\label{subsec:sensitivity}

Finally, to provide a quantitative assessment of the insurer's financial stability in each of the two considered scenarios, let us carry out a sensitivity analysis. For this purpose, we recall that the sensitivity index of a quantity $Q$ which depends differentiably on a parameter $p$ is given by
$$\Upsilon^Q_p = \frac{\partial Q}{\partial p}\cdot \frac{p}{Q},$$
which represents the ratio of the relative change in the value of $Q$ to the relative change in the value of $p$ \cite{ChitnisHymanCushing}.

A sensitivity analysis of the basic reproduction number $\cR_0$ alone, with respect to the epidemiological parameters, is usually carried out in an exact manner, using partial derivatives of $\cR_0$ with respect to the parameters, computed analytically from \eqref{eq:R0}. In this paper, however, we aim to carry out a sensitivity analysis of not only the basic reproduction number $\cR_0$, but also the monthly gross premium $\overline{\cP}^{\text{gross}}$, the minimum loss-preventing start-up capital $\overline{\Gamma}$, and the present value of the end-of-period total profit $\overline{\Pi}_{T/\Delta t}$, with respect to not only the epidemiological parameters but also the insurance-related parameters. Since exact analysis requires explicit expressions of the three financial quantities as functions of the existing parameters, which are not readily available, we shall carry out the sensitivity analysis in a numerical manner, using the following estimate for $\Upsilon^Q_p$:
$$\overline{\Upsilon}^Q_p\left(\Psi\right) = \frac{1}{|\Psi|}\sum_{\psi\in\Psi} \frac{\rho_{Q,p}\left(\psi\right)}{\psi},$$
where $\Psi$ is a set of finitely many signed percentages, and $\rho_{Q,p}\left(\psi\right)$ is the relative change in the value of $Q$ due to a change of $\psi$ in the value of $p$, that is, the post-change value of $Q$ subtracted by the pre-change value of $Q$, divided by the pre-change value of $Q$.

Using $\Psi=\{-10\%,\,-5\%,\,0\%,\,5\%,\,10\%\}$, one obtains the values of $\overline{\Upsilon}^Q_p$ for
\begin{equation}\label{eq:quantities}
Q\in\left\{\cR_0,\,\overline{\cP}^{\text{gross}},\,\overline{\Gamma},\,\overline{\Pi}_{T/\Delta t}\right\}
\end{equation}
and
\begin{equation}\label{eq:parameters}
p\in\left\{\lambda,\,\alpha_1,\,\alpha_2,\,\beta,\,\gamma,\,\mu_1,\,\mu_2,\,i,\,\omega,\,\varphi,\,\cB_H,\,\cB_D,\,\cB_{D^\ast}\right\},
\end{equation}
which are presented in Table \ref{tab:sensdf}. An immediate observation is that for each parameter $p$, the values of $\overline{\Upsilon}^{\cR_0}_p$ in the disease-free and endemic scenarios are equal. In both cases, the basic reproduction number $\cR_0$ depends most sensitively on the natural death coefficient $\mu_1$, with a $1\%$ increase in $\mu_1$ resulting in a $1.00224\%$ decrease in $\cR_0$. 

For the other quantities $Q\in\left\{\overline{\cP}^{\text{gross}},\,\overline{\Gamma},\,\overline{\Pi}_{T/\Delta t}\right\}$ and for each parameter $p$, the values of $\overline{\Upsilon}^{Q}_p$ in the disease-free and endemic scenarios may differ even in their signs. For example, the values of $\overline{\Upsilon}^{\overline{\cP}^{\text{gross}}}_{\lambda}$ is negative in the disease-free scenario and is positive in the endemic scenario, signifying that an increase in $\lambda$ leads to a decrease in $\overline{\cP}^{\text{gross}}$ in the former scenario, and to an increase in $\overline{\cP}^{\text{gross}}$ in the latter scenario. In both scenarios, however, our insurance's monthly gross premium $\overline{\cP}^{\text{gross}}$ depends most sensitively on the disease's incidence coefficient $\beta$, with a $1\%$ increase in $\beta$ leading to a $0.96698\%$ increase in $\overline{\cP}^{\text{gross}}$ in the disease-free scenario, and to a $0.83348\%$ increase in $\overline{\cP}^{\text{gross}}$ in the endemic scenario. On the other hand, in both scenarios, the insurer's minimum loss-preventing start-up capital $\overline{\Gamma}$ depends most sensitively on our insurance's monthly hospitalisation benefit $\mathcal{B}_H$, with a $1\%$ increase in $\mathcal{B}_H$ leading to a $0.67746\%$ increase in $\overline{\Gamma}$ in the disease-free scenario, and to a $0.62895\%$ increase in $\overline{\Gamma}$ in the endemic scenario. Finally, also in both scenarios, the present value of the insurer's end-of-period total profit $\overline{\Pi}_{T/\Delta t}$ depends most sensitively on our insurance's premium surcharge percentage $\varphi$ allocated to the insurer's profit, with a $1\%$ increase in $\varphi$ leading to a $1\%$ increase in $\overline{\Pi}_{T/\Delta t}$ in both scenarios.

\begin{table}\renewcommand{\arraystretch}{1.4}\renewcommand{\tabcolsep}{3pt}
\scalebox{0.925}{    \begin{tabular}{V{4}cV{4}>{\centering\arraybackslash}p{1.5cm}|>{\centering\arraybackslash}p{1.5cm}|>{\centering\arraybackslash}p{1.5cm}|>{\centering\arraybackslash}p{1.5cm}V{4}}\hlineB{4}
        \diagbox{\,$p$\,\,}{\,\,$Q$\,} & $\cR_0$ & $\overline{\cP}^{\text{gross}}$ & $\overline{\Gamma}$ & $\overline{\Pi}_{T/\Delta t}$\\\hlineB{4}
        $\lambda$ & \ar{$1.00000$}	& \ar{$-0.03651$	} &\ar{$0.14368$	} &\arr{$0.25673$}\\\hline
        $\alpha_1$ & \ar{ $0.00000$ } & \ar{ $-0.07325$ } & \ar{ $-0.09691$ } & \arr{ $-0.03878$}\\\hline
        $\alpha_2$ & \ar{ $-0.06866$ } & \ar{ $-0.06904$ } & \ar{	$-0.02963$	} & \arr{ $-0.02716$}\\\hline
        $\beta$ & \ar{ $1.00000$ } & \ar{ $0.96698$ } & \ar{ $0.37118$ } & \arr{ $0.35102$}\\\hline
        $\gamma$ & \ar{ $-0.91092$ } & \ar{ $-0.81215$ } & \ar{ $-0.15733$ } & \arr{ $-0.29950$}\\\hline
        $\mu_1$ & \ar{ $-1.00224$ } & \ar{ $0.19492$} & \ar{$-0.20602$ } & \arr{	$-0.14240$}\\\hline
        $\mu_2$ & \ar{ $-0.02506$ } & \ar{ $-0.35367$ } & \ar{	$0.00549$ } & \arr{	$-0.42841$}\\\hline
        $i$ & \ar{ $0.00000$ } & \ar{ $0.25254$ } & \ar{ $-0.48819$ } & \arr{	$-0.16888$}\\\hline
        $\omega$ & \ar{ $0.00000$ } & \ar{ $0.08696$ } & \ar{ $0.00000$ } & \arr{ $0.00000$}\\\hline
        $\varphi$ & \ar{ $0.00000$ } & \ar{ $0.04348$ } & \ar{ $-0.03734$ } & \arr{ $1.00000$}\\\hline
        $\cB_H$ & \ar{ $0.00000$ } & \ar{ $0.54002$ } & \ar{ $0.67746$ } & \arr{ $0.54002$}\\\hline
        $\cB_D$ & \ar{ $0.00000$ } & \ar{ $0.18694$ } & \ar{ $-0.02174$ } & \arr{ $0.18694$}\\\hline
        $\cB_{D^\ast}$ & \ar{ $0.00000$ } & \ar{	$0.27304$ } & \ar{ $0.34429$ } & \arr{	$0.27304$}\\\hlineB{4}
    \end{tabular}\,\,\,\,\,\,\,\begin{tabular}{V{4}cV{4}>{\centering\arraybackslash}p{1.5cm}|>{\centering\arraybackslash}p{1.5cm}|>{\centering\arraybackslash}p{1.5cm}|>{\centering\arraybackslash}p{1.5cm}V{4}}\hlineB{4}
        \diagbox{\,$p$\,\,}{\,\,$Q$\,} & $\cR_0$ & $\overline{\cP}^{\text{gross}}$ & $\overline{\Gamma}$ & $\overline{\Pi}_{T/\Delta t}$\\\hlineB{4}
        $\lambda$ & \ar{ $1.00000$ } & \ar{ $0.30715$ } & \ar{ $-0.26894$ } & \arr{ $0.37462$}\\\hline
        $\alpha_1$ & \ar{ $0.00000$ } & \ar{ $-0.14336$ } & \ar{ $-0.16772$ } & \arr{ $-0.04043$}\\\hline
        $\alpha_2$ & \ar{ $-0.06866$ } & \ar{	$-0.06395$	} & \ar{ $0.02219$	} & \arr{ $-0.01114$}\\\hline
        $\beta$ & \ar{ $1.00000$ } & \ar{ $0.83348$ } & \ar{ $-0.46500$ } & \arr{ $0.03603$}\\\hline
        $\gamma$ & \ar{ $-0.91092$ } & \ar{ $-0.63554$ } & \ar{ $0.58623$ } & \arr{ $-0.08386$}\\\hline
        $\mu_1$ & \ar{ $-1.00224$ } & \ar{ $-0.08001$ } & \ar{ $0.09159$ } & \arr{  $-0.10917$}\\\hline
        $\mu_2$ & \ar{ $-0.02506$ } & \ar{ $-0.44423$ } & \ar{ $0.26145$ } & \arr{ $-0.52958$}\\\hline
        $i$ & \ar{ $0.00000$ } & \ar{ $0.18283$ } & \ar{ $-0.57299$ } & \arr{ $-0.22825$}\\\hline
        $\omega$ & \ar{ $0.00000$ } & \ar{ $0.08696$ } & \ar{ $0.00000$ } & \arr{ $0.00000$}\\\hline
        $\varphi$ & \ar{ $0.00000$ } & \ar{ $0.04348$ } & \ar{ $-0.06751$ } & \arr{ $1.00000$}\\\hline
        $\cB_H$ & \ar{ $0.00000$	} & \ar{ $0.62895$} & \ar{	$0.67498$ } & \arr{ $0.62895$ }\\\hline
        $\cB_D$ & \ar{ $0.00000$ } & \ar{ $0.05290$ } & \ar{ $-0.01894$	} & \arr{ $0.05290$}\\\hline
        $\cB_{D^\ast}$ & \ar{ $0.00000$} & \ar{	$0.31815$	} & \ar{ $0.34396$	} & \arr{ $0.31815$}\\\hlineB{4}
    \end{tabular}}\smallskip

\caption{\label{tab:sensdf} Estimates $\overline{\Upsilon}^Q_p\left(\{-10\%,\,-5\%,\,0\%,\,5\%,\,10\%\}\right)$ for the sensitivity indices $\Upsilon^Q_p$ of $Q$ with respect to $p$, for all $Q$ and $p$ given by \eqref{eq:quantities} and \eqref{eq:parameters}, in our disease-free scenario $\cR_0<1$ (left) and in our endemic scenario (right).}
\end{table}

\section{Conclusions and future research}\label{sec:conclusions}

Continuing the existing studies on epidemic-model-based health insurances \cite{Feng,FengGarrido,NkekiEkhaguere,Hainaut,AtatalabEtAl,Nam,Guerra,ZhaiEtAl}, we have presented a design and financial analysis of a health insurance based on an SIH-type epidemic model. Specifically, we have, firstly, constructed the model in continuous form, analysed the model dynamically, and formulated the financial quantities involved in the insurance. Subsequently, we have constructed a discrete form of the model using the forward Euler method, carried out the analogous dynamical analysis for this discrete model, and developed discrete estimates for the previously formulated financial quantities. Using the same discrete model, we have simulated two numerical scenarios, disease-free and endemic, corresponding to two different sets of parameter values derived from a range of sources, differing only in the value of the disease's incidence coefficient. The results indicated that, in the endemic scenario, our insurance's gross premium is more than three times higher, the insurer's minimum loss-preventing start-up capital is lower, and the insurer's end-of-period total profit is higher, compared to the respective values in the disease-free scenario. Finally, we have conducted a sensitivity analysis, which revealed that, in both scenarios, the disease's basic reproduction number, the monthly gross premium, the minimum start-up capital, and the total profit depend most sensitively on, respectively, the population's natural death coefficient, the disease's incidence coefficient, the hospitalisation benefit, and the premium surcharge percentage allocated to profit.

Our work is open to further advancements. First, one could investigate whether the global stability of both equilibria in our continuous and discrete models could be established using, e.g., Lyapunov functions. Next, with regards to our models' construction, various paths for extension are possible. For instance, one could replace our models' bilinear form of incidence rate and linear form of treatment rate with more realistic forms, such as the Holling and Beddington-DeAngelis forms \cite{Nilam}. In addition, to capture the epidemiological and actuarial impacts of the individuals' heterogeneous ages, one could construct age-structured versions of our models, dividing each compartment into subcompartments representing different age groups \cite{RamSchaposnik,CuevasMaraverEtAl,CanabarroEtAl,MistryEtAl,TaboeEtAl}. Furthermore, the case in which only a portion of the population is insured may also be of interest. Finally, more realistic formulations of the insurance-related quantities could be obtained by incorporating, e.g., taxes, different hospitalisation benefit amounts and policies for different treatment types, and compensations for individuals whose contracts are terminated at the end of the insurance's availability period, and by assuming non-constant interest rates described by stochastic models, such as the Vasicek and Cox-Ingersoll-Ross models \cite{Vasicek,CoxEtAl}.

\section*{Acknowledgements}

The initiative to write this article stems from the participation of the last three authors in the Actuarial Case Competition 2024 at the Institute of Technology Bandung, Indonesia. The authors thus thank the competition's organisers and contributors, especially the judges Jasen Theodorus and Dila Puspita for their valuable suggestions.

%During the writing of this article, the authors used ChatGPT in order to improve language and readability. After using this tool, the authors reviewed and edited the content as needed and take full responsibility for the content of the publication.

\section*{Disclosure statement}

No potential conflict of interest was reported by the authors.

%This article is a developed version of the report written by the first three authors, under the supervision of the last two authors, as a team participating in the Actuarial Case Competition 2024, an event hosted by the Institute of Technology Bandung (ITB). The authors thus thank the organisers of the event.

%%%%%%%%%%%%%%%%%%%%%%%%%%%%%%%%%%%%%%%%%%%%%%%%%%%%%%%%%%%%%%%%%%%%%%%%%%%

\end{document}